\documentclass{article}
\usepackage{authblk}

\usepackage[T2A]{fontenc}
\usepackage[english]{babel}
\usepackage{amssymb}
\usepackage{multirow}
\usepackage{fancyhdr}

\large

\newtheorem{Prop}{Proposition}
\newtheorem{Theo}{Theorem}

\title{\textbf{Generalizations of Mezhirov's game semantics to predicate superintuitionistic logics and the Casari formula}}
\author{Ivan O. Pyltsyn}
\date{}
\affil{HSE University}

\begin{document}
	\maketitle
	\begin{abstract}
		\noindent Game semantics allows us to look at basic logical concepts from another side. This approach to logic has a long history, there are plenty of different types of games: provability games, semantic games, etc [10,11]. And there is an interesting type of provability games called Mezhirov's game proposed by Iliya Mezhirov for intuitionistic logic of propositions ($IPC$) and Grzegorczyk modal logic ($Grz$) [1,2]. This idea was developed in many different directions; for example, in 2008 in the joint paper with  N. Vereschagin a game semantics was given for affine and linear logic [3]. Independently G. Japaridze worked on game semantics for linear logic [4]. Mezhirov's games for minimal propositional logic ($MPC$), logic of functional frames ($KD!$) and logic of serial frames ($KD$) were introduced in 2021 by A. Pavlova [5].
	\end{abstract}
	\setcounter{page}{1}
	\section{Introduction}
	\quad Mezhirov's game semantics for intuitionistic logic is interesting because of its simplicity and strong connection with Kripke semantics and Kripke models. The game between Opponent and Proponent starts with a formula $\varphi$. And Proponent has a winning strategy iff $\varphi$ is an intuitionistic tautology. The connection between the game and Kripke models manifests itself in building strategy for Opponent from a Kripke model (Opponent "walks" from one world of a model to another) and in the reconstruction of a model from Opponent's winning strategy (in which there exists a world where $\varphi$ is false). And these procedures are connected to each other.
	
	In my study, I try to generalize Mezhirov's result in two directions: to generalize to intuitionistic logic of predicates (introduce a game between Opponent and Proponent with at least the same connection with Kripke models or with special classes of them) and to the case of a connection not only between the game and tautologies of logic ($\vDash\varphi$), but also between the game and entailment from infinite sets of formulas ($T\vDash\varphi$).
	
	The purpose of building such game was to describe a game (winning strategies in it) by a semantic consequence defined by some class of predicate Kripke frames ($\varphi$ is a semantic consequence of $T$ in some class $C$ of predicate Kripke frames iff for each Kripke model, based on a frame from the class $C$, if all formulas from $T$ is true everywhere in this model, then $\varphi$ is true in each world of the model), and additionally get a description of this semantic consequence by winning strategies in the game. I initially thought about just logic of all Kripke models, i.e. it would be a game for intuitionistic logic of predicates directly. But it turned out that in such case some fundamental problems arise and it is natural to expand the logic (to use a smaller class of Kripke frames). Moreover, description of such variations (not just logic of all Kripke models) could be useful, since, in general, Kripke semantics for superintuitionistic predicate logic is rather weak (e.g. [9]). And I managed to get a description (based on the game I built) for several variations.
	
	\section{Main definitions}
	\quad Let $\Omega$ be the elementary intuitionistic language (without function symbols; the set of logical connectives will be $\{\to, \wedge, \lor, \bot\}$, where $\neg A$ will be considered as $A\to\bot$), and we will use Kripke models for intuitionistic logic of predicates [6,7] (I will call sets of constants in each world "individual domains" (or "the set of objects") and use symbol $\Delta$). For the set of formulas $\Gamma$ and set of objects (constants) $\Delta$ let $\mathcal{F}(\Gamma, \Delta)=\{P[c_{1}, ..., c_{n}] | P[x_{1}, ..., x_{n}]$ is a subformula of some formula from $\Gamma$ and free variables of it are only $x_{1}, ..., x_{n}; c_{i}\in\Delta\}$ (so $\mathcal{F}$ in some ways is a set of all "subformulas" of formulas from $\Gamma$). Players	Opponent and Proponent will be associated with their sets $\mathcal{O}$ and $\mathcal{P}$. The position in the game is a triple $\mathcal{C}=(\mathcal{O}, \mathcal{P}, \Delta)$. In each position $\mathcal{C}$: $\mathcal{O}$ and $\mathcal{P}$ are subsets of $\mathcal{F}(\Gamma, \Delta)$, where $\Delta$ is taken from $\mathcal{C}$ (and can only expand in the game process) and $\Gamma$ is fixed at the beginning of the game and does not change until the end and equals to $\mathcal{O}_{0}\cup\{\varphi\}$ (where $\mathcal{C}_{0}=(\mathcal{O}_{0}, \{\varphi\}, \Delta_{0})$ is a starting position; $\Delta_{0}$ is an exact set of all constants contained in formulas from $\Gamma$). Proponent moves by adding new formulas from $\mathcal{F}$ to $\mathcal{P}$, Opponent moves by expanding $\Delta$ (he can add nothing to $\Delta$ if he wants; and he can add to $\Delta$ not just constants from $\Omega$, but any new elements) and then adding new formulas from $\mathcal{F}$ to $\mathcal{O}$.
	
	The only thing left to define is who must move in a position $\mathcal{C}$. To do that, let us firstly define the notion of truth relation $\Vdash$ in $\mathcal{C}$ for formulas from $\mathcal{F}(\Gamma, \Delta)$:
	
	\begin{table}[h!]
		\centering
		\begin{tabular}{l}
			$\mathcal{C}\nVdash\bot$ \\
			$\mathcal{C}\Vdash A[c_{1}, ..., c_{n}] \rightleftharpoons A[c_{1}, ..., c_{n}]\in\mathcal{O} $ \\
			$\mathcal{C}\Vdash\varphi\star\psi\rightleftharpoons\varphi\star\psi\in\mathcal{O}\cup\mathcal{P}$ and $ (\mathcal{C}\Vdash\varphi)\star(\mathcal{C}\Vdash\psi)$, $\star\in \{\to, \wedge, \lor\}$ \\
			$\mathcal{C}\Vdash qxP[x]\rightleftharpoons qxP[x]\in\mathcal{O}\cup\mathcal{P}$ and $q\alpha\in\Delta$($\mathcal{C}\Vdash P[\alpha]$), $q\in\{\exists, \forall\}$
		\end{tabular}
	\end{table}
	
	\noindent where $A$ is a predicate symbol, $arity(A)=n$, $c_{i}\in\Delta$, $P$ - formula with only one free variable. A star in the case of $(\mathcal{C}\Vdash\varphi)\star(\mathcal{C}\Vdash\psi)$ means logical meta connective and behaves like a classical connective (the same holds for $q$ in $q\alpha\in\Delta$).
	
	Let us call formulas from $\mathcal{O}\cup\mathcal{P}$ marked formulas (only they can be true in the current position).
	
	Let us call a formula from $\mathcal{P}$ Proponent's mistake if it is false in the current position (the same for $\mathcal{O}$ and Opponent). If Opponent has no mistakes but Proponent has, then Proponent moves. Otherwise, Opponent must move. And if after a turn of a fixed player this player must move again, he loses and the game ends. If the game goes on infinitely (each player manages to pass a turn to the other player each turn), Proponent wins (more formally, the game is a possibly infinite sequence of positions, where adjacent positions behaves as some player in his turn made his move as described from a position to the next one).
	
	\section{Classes of frames}
	\quad In the final analysis, the mentioned game could describe to varying degrees several classes of frames. Let us give them abbreviations. 
	
	\noindent $N\leftrightharpoons$ class of all Noetherian Kripke frames (frames with no strictly increasing infinite sequences of worlds).
	
	\noindent $Cas\leftrightharpoons$ class of all Kripke frames in which Casari's formula $\forall x[(P[x]\to \forall xP[x])\to \forall xP[x]]\to\forall xP[x]$ is valid; let us call it Casari's class (Kripke frame is from Casari's class iff in every countable sequence of worlds $\omega_{i}$ their individual domains $\Delta_{i}$ remain finite and stabilize; here I consider Kripke frame to be not just partially ordered set, partially ordered set with individual domain in each world (which is usually a part of a model)) [8].
	
	\noindent $Fin\leftrightharpoons$ class of all finite Kripke frames.
	
	If $X$ is a class of Kripke frames, then let $finX$ be a subclass of $X$ of all frames from $X$ with only finite individual domains in each world.
	
	In particular, $Cas\supset N\supset Fin$, $Cas\supset finCas$, etc.
	
	And let an abbreviation for semantic consequence given by a class of frames $X$ be $\vDash_{X}$.
	
	\section{Examples}
	\quad It seems to me that, informally, this game (and Mezhirov's game for propositional intuitionistic logic) could be understood as follows: Opponent is trying to build a theory that belies Proponent's assertion that $\phi$ follows from $\mathcal{O}_{0}$ (or, in the case of $\mathcal{O}_{0}=\varnothing$, is trying to build a theory that shows that Proponent's thesis ($\varphi$) is not valid in general). And this theory must be coherent (Opponent must have no mistakes), otherwise his approach is considered unsuccessful.  
	
	Now let us consider several examples of the game.
	
	In the first game $\mathcal{C}_{0}=(\varnothing, \{\varphi\}, \varnothing)$, where $\varphi=\forall y\exists x(P[x]\to P[y])$. Because $\Delta$ is empty, there are no formulas in $\mathcal{F}$ of the form $\exists x(P[x]\to P[c])$, so Proponent has no mistakes, it is Opponent's turn. It is enough for him to just expand $\Delta$, and it will be Proponent's turn. Proponent takes all formulas of the kind $\exists x(P[x]\to P[c])$ and $P[c]\to P[c]$ and passes turn to Opponent. He will do the same (expand $\Delta$) and the game goes on infinitely. Informally, Opponent each turn asks Proponent, what he thinks about this new object. And Proponent shows that even for this new object he can prove his thesis. This intuitive view, I suppose, shows that the universal quantifier here is closely related to the intuitionistic universal quantifier (but not completely the same).
	
	In the second game $\mathcal{C}_{0}=(\varnothing, \{\varphi\}, \{c\})$, where $\varphi=\neg P[c]\to\neg\exists xP[x]$. $\varphi$ is an implication, both antecedent and consequent of it is not marked, therefore false in the current position. So $\varphi$ is true, it is Opponent's turn. He expand $\Delta$ to $\{c, \alpha\}$ and add to $\mathcal{O}$ formulas $\neg P[c]$, $\exists xP[x]$, $P[\alpha]$. He might not add $\exists xP[x]$ to $\mathcal{O}$ and turn would still be passed to Proponent. But in this case Proponent would have an opportunity to add to $\mathcal{P}$ $\neg \exists xP[x]$ and make this formula true in position (because $\exists xP[x]$ would not be marked), and Opponent still would have needed to add $\exists xP[x]$. After that, Proponent will not be able to pass the turn to Opponent, therefore, he will lose. 
	
	In the third game let $\mathcal{C}_{0}=(\varnothing, \{\varphi\}, \varnothing)$, where $\varphi=\forall x[(P[x]\to \forall xP[x])\to \forall xP[x]]\to\forall xP[x]$ (\emph{Casari's schema or Casari's formula}). Again $\varphi$ is an implication, it is Opponent's turn. He needs to make sending false, so he expand $\Delta$ and add to $\mathcal{O}$ all formulas $(P[\alpha]\to \forall xP[x])\to \forall xP[x]$ and antecedent of the $\varphi$: $\forall x[(P[x]\to \forall xP[x])\to \forall xP[x]]$. Then Propopent creates mistakes for Opponent by adding to $\mathcal{P}$ all formulas $(P[\alpha]\to \forall xP[x])$. To get rid of mistakes, Opponent needs to add all $P[\alpha]$, and then Proponent just add to $\mathcal{P}$ $\forall xP[x]$. The only thing Opponent can do now is to expand $\Delta$ and repeat everything again. As we can see, this is the winning strategy for Proponent, but $\varphi$ is not true in all Kripke models. This is why we are interested in $Cas$.
	
	\section{Relation between the game and $N$}
	
	\begin{Theo}
		Proponent has a winning strategy in position $\mathcal{C}_{0}=(\mathcal{O}_{0}, \{\varphi\}, \Delta_{0})$ iff $\mathcal{O}_{0}\vDash_{N}\varphi$.
	\end{Theo}
	
	\noindent Here $\mathcal{O}_{0}$ - arbitrary subset of a set of all closed formulas of language $\Omega$.
	
	We will prove this theorem in two stages: firstly, we will show that \emph{$\mathcal{O}_{0}\nvDash_{N}\varphi\Rightarrow$ starting from position $\mathcal{C}_{0}$ Opponent wins}, having built a strategy for Opponent from a countermodel (due to which we now that $\mathcal{O}_{0}\nvDash_{N}\varphi$); secondly: \emph{$\mathcal{O}_{0}\vDash_{N}\varphi\Rightarrow$ starting from position $\mathcal{C}_{0}$ Proponent wins}, because he can pass turn to Opponent each time by just adding to $\mathcal{P}$ all formulas from $\mathcal{F}$ that follows from $\mathcal{O}$ (we will show that this strategy for Proponent works in such case).
	
	\subsection{Strategy for Opponent}
	
	\begin{Prop}
		Let the current position be $\mathcal{C}=(\mathcal{O}, \mathcal{P}, \Delta)$, and in a model $M$ based on a Kripke frame $K$ from the class $N$ (model and frame over language $\Omega$ extended by all objects from $\Delta$ that is not constants from $\Omega$) in a world $\omega$ all formulas from $\mathcal{O}$ is true, and some formula from $\mathcal{P}$ is not. Then from this position Opponent has a winning strategy.
	\end{Prop}
	
	\noindent \emph{Proof.} Firstly, we can consider the model (and frame) to be a cone of the world $\omega$ (cone of Noetherian frame is still Noetherian frame). Hence, in the whole model all $\mathcal{O}$ is true.
	
	Let it be Opponent's turn in the current position. Choose any maximal world from our model among those in which something from $\mathcal{P}$ is falsified (this set is not empty because in $\omega$ something from $\mathcal{P}$ is falsified and it is Noetherian Kripke frame, there is no infinite sequences). Let us denote this world as $w$. Now Opponent as his move extends $\Delta$ to $\Delta'$ by adding all elements from $\Delta_{w}$ that is not an interpretation of some element from $\Delta$ (by interpretation I mean interpretation of constants of our language in the model) and add to $\mathcal{O}$ all formulas from $\mathcal{F}(\Gamma, \Delta')$ that is true in $w$. Let us show that in a resulting position $\mathcal{C}'=(\mathcal{O}', \mathcal{P}, \Delta')$ we have $\mathcal{C}'\Vdash\varphi\Leftrightarrow w\Vdash\varphi$ (for formulas $\varphi\in\mathcal{F}$; I will not mention this phrase each time later). For unmarked formulas it is true because by definition of $\Vdash$ in position we have $\mathcal{C}'\nVdash\varphi$, and since they are unmarked, they are not in $\mathcal{O}'$, which is the set of all true formulas in $w$, so they are not true in $w$. Now use induction on the constraction of $\varphi$. If it is an atomic formula, then $\mathcal{C}'\Vdash\varphi$ is equivalent to $\varphi\in\mathcal{O}'$, which is the set of all true formulas in $w$. If $\varphi$ is not atomic and unmarked, we already now that proposition holds. So let us consider that $\varphi$ is not atomic and marked. Then everywhere above $w$ this formula is true (because $w$ is a maximal world, in which something from $\mathcal{P}$ is falsified, and $\mathcal{O}'$ is true in $w$, so it is true everywhere above), hence validity of $\varphi$ is determined by validity of its subformulas classically (even in case of $\varphi=A\to B$ and $\varphi=\forall xP[x]$), and in the same way truth defines in the position for marked formulas according to its subformulas. Proposition holds for subformulas by induction hypothesis, therefore for such $\varphi$.    
	
	Then we have $\mathcal{C}'\Vdash\varphi\Leftrightarrow w\Vdash\varphi$. And we know that all $\mathcal{O}'$ is true in $w$ and something from $\mathcal{P}$ is not. Hence, in $\mathcal{C}'$ Opponent has no mistakes but Proponent has. It is Proponent's turn. 
	
	We actually showed that in the situation described in proposition Opponent can pass a turn to Proponent by "moving" from the world $\omega$ to the world $w$. When Proponent passes turn back to Opponent, he will do the same (he can use the same model (or just cone of the world $w$) for finding new world according to which he will make his move, because $\mathcal{P}$ can only be extended) and still be able to pass turn to Proponent. Moreover, each Opponent's turn will be associated with some world of the model from which we started. This model is based on the Noetherian Kripke frame, hence it is not possible for the game to be infinite (since each new turn of the Opponent gives us a new world above the previous one). As was mentioned, Opponent can always pass a turn. Therefore, the last turn was Proponent's turn, he did not manage to pass turn to Opponent (after his turn, it is his turn again), which means that he lost. Opponent won.
	
	\begin{flushright}
		$\square$
	\end{flushright}  
	
	Now we can apply proven proposition to the starting position $\mathcal{C}_{0}$ and from the model which shows that $\mathcal{O}_{0}\nvDash_{N}\varphi$ get winning strategy for Opponent.
	
	\subsection{Justification of the Proponent's strategy}
	\begin{Prop}
		Let the current position be $\mathcal{C}=(\mathcal{O}, \mathcal{P}, \Delta)$ and $\mathcal{O}\vDash_{N}\mathcal{P}$ ($\forall\varphi\in\mathcal{P}: \mathcal{O}\vDash_{N}\varphi$; semantic consequence given by a models over language $\Omega$ extended by all objects from $\Delta$ that is not constants from $\Omega$). Then Proponent has a winning strategy from this position.
	\end{Prop}
	
	\noindent \emph{Proof.} It is enough for Proponent to always have an opportunity to pass a turn to Opponent. Hence, without limiting the generality, let us consider that in the current position it is Proponent's turn (because $\mathcal{O}$ can only expand, and condition $\mathcal{O}\vDash_{N}\mathcal{P}$ will still hold). Then Proponent extend $\mathcal{P}$ to $\mathcal{P}'$ by adding all formulas that semantically follows from $\mathcal{O}$. We need to show that in the new position $\mathcal{C}'=(\mathcal{O}, \mathcal{P}', \Delta)$ it will be Opponent's turn.
	
	Each formula $\varphi_{i}\in\mathcal{F}\setminus\mathcal{P}'$ does not follow fron $\mathcal{O}$, hence exist a model $M_{i}$ over Noetherian Kripke frame $K_{i}$ and world $\omega_{i}$ (and we can consider this model to be a cone of the world $\omega_{i}$) in which all $\mathcal{O}$ is true (therefore in the whole model $M_{i}$) and $\omega_{i}\nVdash\varphi_{i}$. In these models interpretations of different elements from $\Delta$ could be the same. In such case we should modify $M_{i}$ by adding to individual domains of each world new elements for these elements to make interpretations unique (truth for the atomic formulas with new elements the same as truth for atomic formulas where we replaced all new elements by the old interpretations of elements from $\Delta$ for which we added these new elements). And also by set-theortic simple managment we can consider that interpretation of each element from $\Delta$ is this element (so $\Delta$ is a subset of any individual domain). Let us build a new model: add to a disjoint union of all $M_{i}$ root $\omega$ with individual domain $\Delta$, and truth for the atomic formulas: everywhere above $\omega$ truth is the same as in $M_{i}$, and in $\omega$ atomic formula is true iff it belongs to $\mathcal{O}$. Monotonicity of individual domains and truth is preserved ($\mathcal{O}$ is true in each $M_{i}$), and because each $K_{i}$ is Noetherian, this model is based on Noetherian Kripke frame. 
	
	Let us show that $\mathcal{C}'\Vdash\varphi\Leftrightarrow\omega\Vdash\varphi$. For unmarked formulas (they are not in $\mathcal{P}'$, so unmarked formulas are exactly all $\varphi_{i}\in\mathcal{F}\setminus\mathcal{P}'$) proposition holds because they by definition are not true in the position and falsified somewhere above $\omega$, hence they are not true in $\omega$. Now again use induction on the construction of the formula $\varphi$. If it is atomic, the truth of this formula in the position is equivalent to belonging to $\mathcal{O}$ and that is how we defined the truth for atomic formulas in $\omega$. Now let us consider $\varphi$ to not be atomic and marked. This formula is marked (which means that $\varphi\in\mathcal{O}\cup\mathcal{P}'$) then it is true everywhere above $\omega$, hence its validity is determined by its subformulas classically, as truth in the position. And, by induction hypothesis, the truth of subformulas is the same in position and in $\omega$. Therefore, the same holds for $\varphi$.
	
	Now we have $\mathcal{C}'\Vdash\varphi\Leftrightarrow\omega\Vdash\varphi$. If in $\omega$ not all formulas from $\mathcal{O}$ are true, then the same holds for the current position, which means that Opponent has mistakes, it is his turn. If all $\mathcal{O}$ is true in the $\omega$, than the same holds for the whole model, and because this model is based on Noetherian Kripke frame, all $\mathcal{P}'$ is also true in the model, hence it is true in $\omega$, so Proponent has no mistakes in the current position. It is Opponent's turn.
	
	We've managed to show that in such case Proponent can pass the turn to Opponent. Moreover, condition $\mathcal{O}\vDash_{N}\mathcal{P}'$ will still hold, even after Opponent's turn (he possibly will extend $\mathcal{O}$ and $\Delta$). It means that Proponent can repeat this strategy each time. He will always be able to pass a turn to Opponent. Hence, he wins.
	
	\begin{flushright}
		$\square$
	\end{flushright}  
	
	Applying this proposition to the starting position $\mathcal{C}_{0}$ we get what we wanted.   
	
	\section{Relation between game and $Cas$}
	\begin{Theo}
		For finite $\mathcal{O}_{0}$: Proponent has a winning strategy in position $\mathcal{C}_{0}=(\mathcal{O}_{0}, \{\varphi\}, \Delta_{0})$ iff $\mathcal{O}_{0}\vDash_{Cas}\varphi$.
	\end{Theo}
	
	In particular, this implies $T\vDash_{Cas}\varphi\Leftrightarrow T\vDash_{N}\varphi$ for finite $T$.
	
	Analogously, let us prove this theorem in two stages. And one part of it we can get from the previous theorem: \emph{$\mathcal{O}_{0}\vDash_{Cas}\varphi\Rightarrow\mathcal{O}_{0}\vDash_{N}\varphi\Rightarrow$ starting from position $\mathcal{C}_{0}$ Proponent wins}. We only need to prove another part: \emph{$\mathcal{O}_{0}\nvDash_{Cas}\varphi\Rightarrow$ starting from position $\mathcal{C}_{0}$ Opponent wins}.
	
	\subsection{Strategy for Opponent}
	\begin{Prop}
		Let the current position be $\mathcal{C}=(\mathcal{O}, \mathcal{P}, \Delta)$ ($\Gamma$ is finite in this game), and in a model $M$ based on a Kripke frame $K$ from the class $Cas$ (model and frame over language $\Omega$ extended by all objects from $\Delta$ that is not constants from $\Omega$) in a world $\omega$ all formulas from $\mathcal{O}$ is true, and some formula from $\mathcal{P}$ is not. Then from this position Opponent has a winning strategy.
	\end{Prop}
	
	\noindent \emph{Proof.} Again, we can consider $M$ (and $K$) to be a cone of the world $\omega$ (frame will still be Casari's frame). And in this case all formulas from $\mathcal{O}$ are true in $M$.
	
	Let it be Opponent's turn in the current position. We will build a strictly increasing sequence of worlds (informally, this procedure can be seen as Opponent is planning his move). $\omega_{0}=\omega$. Having $\omega_{n}$: from the worlds above $\omega_{n}$ we choose one in which something from $\mathcal{P}$ is falsified and its individual domain is strictly extends individual domain of the world $\omega_{n}$. If there is no such worlds, we will choose among the worlds above $\omega_{n}$ one in which something from $\mathcal{P}$ is falsified (hence, this world has the same individual domain $\Delta_{\omega_{n}}$) and truth for the formulas from $\mathcal{F}(\Gamma, \Delta_{\omega_{n}})$ is different from the world $\omega_{n}$. If there is no such worlds, we stop our sequence. Let us consider that builded sequence is infinite. This is strictly increasing infinite sequence of worlds of Casari's frame, hence exists a world $\omega_{n}$ from which individual domains in sequence stopped changing and remain finite. This means that from that step (by construction of the sequence) truth for the formulas from $\mathcal{F}(\Gamma, \Delta_{\omega_{n}})$ must change each time in the sequence, but because of monotonicity of truth and finiteness of $\mathcal{F}(\Gamma, \Delta_{\omega_{n}})$ (because $\Gamma$ and $\Delta_{\omega_{n}}$ is finite), truth cannot change infinitely many times, it is a contradiction. 
	
	Hence, this sequence is finite. Let the last world of this sequence be $\omega_{m}=w$. In $w$ something from $\mathcal{P}$ is falsified and above this world (by the construction of the sequence (we know, that on the last step we stopped building the sequence)) only worlds of two types: $X=\{\omega\geqslant w | \omega\Vdash\mathcal{P}\}$ and $Y=\{\omega\geqslant w | \Delta_{\omega}=\Delta_{w}, \forall\psi\in\mathcal{F}(\Gamma, \Delta_{w})[w\Vdash\psi\Leftrightarrow\omega\Vdash\psi]\}$ ($X$ is a set of worlds of the first type, and $Y$ - set of worlds of the second type).
	
	Now Opponent similarly to the Noetherian case as his turn extends $\Delta$ to $\Delta'$ by adding all elements from $\Delta_{w}$ that is not an interpretation of some element from $\Delta$ and add to $\mathcal{O}$ all formulas from $\mathcal{F}(\Gamma, \Delta')$ that is true in $w$. Let us show that in new position $\mathcal{C}'=(\mathcal{O}', \mathcal{P}, \Delta')$ holds: $\mathcal{C}'\Vdash\varphi\Leftrightarrow w\Vdash\varphi$. For unmarked formulas proposition is true because $\mathcal{C}'\nVdash\varphi$ and since they are unmarked, they are not in $\mathcal{O}'$, which is exect set of all true formulas in $w$. Now use induction on the constraction of $\varphi$. If $\varphi$ is an atom, then its validity in the position is equivalent to belonging to $\mathcal{O}'$ - the exact set of truth formulas in $w$. Let us consider $\varphi$ to be not atomic and marked. $\varphi$ is marked, hence it is true in all worlds from $X$ ($\mathcal{O}'$ is true everywhere above $w$) and truth in worlds from $Y$ is the same as in $w$. Therefore, the validity of $\varphi$ in $w$ is determined by the validity of its subformulas classically, and in the same way truth defines in the position for marked formulas according to its subformulas. Proposition holds for subformulas by induction hypothesis, hence for such $\varphi$.
	
	Therefore we have $\mathcal{C}'\Vdash\varphi\Leftrightarrow w\Vdash\varphi$. And we know that $w\Vdash\mathcal{O}'$ and something from $\mathcal{P}$ is falsified in $w$. Hence, it is Proponent's turn. 
	
	We actually showed that in the situation described in the proposition Opponent can pass the turn to Proponent by "moving" from the world $\omega$ to the world $w$ (analogously to the Noetherian case). When Proponent passes turn back to Opponent, he will do the same (he can use the same model (or just cone of the world $w$) for finding new world according to which he will make his move, because $\mathcal{P}$ can only be extended) and still be able to pass turn to Proponent. Moreover, each Opponent's turn will be associated with some world of the model from which we started. This model is based on Casari's Kripke frame, hence it is not possible for the game to be infinite: each new turn of the Opponent gives us a new world above the previous one, therefore we will get an infinite strictly increasing sequence of worlds. Then, after some turn Opponent will stop extending $\Delta$, individual domains $\Delta_{w}$ in sequence will also stop changing and remain finite and from some step, since $\mathcal{F}(\Gamma, \Delta_{w})$ is finite, truth will stop changing, which means (because each world corresponds to Opponent's turn) that, from some turn, Opponent will not do anything and still will pass turn to Proponent. This is the same as Proponent did not manage to pass a turn to Opponent, which means that he lost, and the game is not infinite. And since Opponent always has an opportunity to pass a turn to Proponent, he will win.	
	\begin{flushright}
		$\square$
	\end{flushright}  
	
	Now, again, we can apply proven proposition to the starting position $\mathcal{C}_{0}$ and from the model which shows that $\mathcal{O}_{0}\nvDash_{Cas}\varphi$ get winning strategy for Opponent.
	
	\section{Relation between the game and classes $fin\dots$}
	\quad Now let us take a look at the same game, but we will consider only finite $\Delta_{0}$ and Opponent will be allowed to extend $\Delta$ by adding only finitely many elements in a single turn (let us call this "finite variation of the game"). By analogous reasoning, we can obtain following results for $finN$ and $finCas$:
	
	\begin{Theo}
		In the finite variation, Proponent has a winning strategy in position $\mathcal{C}_{0}=(\mathcal{O}_{0}, \{\varphi\}, \Delta_{0})$ (with possibly infinite $\mathcal{O}_{0}$, but with only finite $\Delta_{0}$) iff $\mathcal{O}_{0}\vDash_{finN}\varphi$.
	\end{Theo}
	\begin{Theo}
		In the finite variation, Proponent has a winning strategy in position $\mathcal{C}_{0}=(\mathcal{O}_{0}, \{\varphi\}, \Delta_{0})$ (with only finite $\mathcal{O}_{0}$) iff $\mathcal{O}_{0}\vDash_{finCas}\varphi$.
	\end{Theo}
	
	\noindent Moreover, we can obtain following result (proof is similar to proof for finN):
	
	\begin{Theo}
		In the finite variation, Proponent has a winning strategy in position $\mathcal{C}_{0}=(\mathcal{O}_{0}, \{\varphi\}, \Delta_{0})$ (with only finite $\mathcal{O}_{0}$) iff $\mathcal{O}_{0}\vDash_{finFin}\varphi$.
	\end{Theo}
	
	\noindent In particular: $T\vDash_{finCas}\varphi\Leftrightarrow T\vDash_{finN}\varphi\Leftrightarrow T\vDash_{finFin}\varphi$ for finite $T$.
	
	I should also mention that from the winning strategy for Opponent in the finite variation of the game that starts from position $\mathcal{C}_{0}=(\mathcal{O}_{0}, \{\varphi\}, \Delta_{0})$ we can construct a countermodel in which all $\mathcal{O}$ is true, and $\varphi$ is falsified. This procedure is strongly connected with building strategy for Opponent from a countermodel. I will not describe it here since it is mostly similar to the analogous procedure in [1].
	
	\section{Conclusion and future work}
	\quad As I mentioned, the main goal of this study was to find a game with strong connection with Kripke models. Partially, this has been achieved; in addition, some connections have been established between weak entailment of logics of some classes. I suppose, the next step that would be to find a triple: a class of Kripke frames, a game semantics and a calculus (probably, an infinitary sequent calculus) with the same strong entailment (entailment from not only finite, but from any sets of formulas). In this case, it is better to take a simpler class of Kripke frames in terms of the possible calculus for this class. So, because of this, Casari's class looks better than the Noetherian class. Therefore, it is natural to try to change rules of the game to get the same strong entailment as in logic of all Kripke frames from Casari's class.
	
	\section*{Acknowledgements}
	\noindent This research was done under the supervision of Associate Professor Daniyar Shamkanov.
	
	\noindent The work was partially supported by the HSE Academic Fund Programme (Project: Research and Study Group number 23-00-022).

\end{document}